\documentclass[reqno,12pt]{amsart}
\usepackage{amsmath, latexsym, amsfonts, amssymb, amsthm, amscd}

\numberwithin{equation}{section}

\newtheorem{theorem}{Theorem}[section]
\newtheorem{lemma}[theorem]{Lemma}

\newtheorem{corollary}[theorem]{Corollary}
\newtheorem{rem}[theorem]{Remark}
\newtheorem{definition}[theorem]{Definition}

\DeclareMathSymbol{\leqslant}{\mathalpha}{AMSa}{"36} 
\DeclareMathSymbol{\geqslant}{\mathalpha}{AMSa}{"3E} 
\DeclareMathSymbol{\eset}{\mathalpha}{AMSb}{"3F}     
\renewcommand{\leq}{\;\leqslant\;}                   
\renewcommand{\geq}{\;\geqslant\;}                   

\newcommand{\R}{\mathbb{R}}
\newcommand{\Z}{\mathbb{Z}}
\newcommand{\N}{\mathbb{N}}

\title[Strong localization and macroscopic atoms for directed polymers]{Strong localization and  macroscopic atoms for directed polymers}

\author{Vincent Vargas}
\thanks{Partially
   supported by CNRS (UMR 7599
``Probabilit{\'e}s et Mod{\`e}les
Al{\'e}atoires'')}
        
\begin{document}

\maketitle
\begin{center}
{\footnotesize \noindent
 Universit{\'e} Paris 7,\\
Math{\'e}matiques, case 7012,\\ 2, place Jussieu, 75251 Paris, France}

{\footnotesize \noindent e-mail: \texttt{vargas@math.jussieu.fr}}
\end{center}

\begin{abstract}
In this article, we derive strong localization results for directed polymers
in random environment. We show that at "low temperature" the polymer measure
is asymptotically concentrated at a few points of macroscopic mass
(we call these points $\epsilon$-atoms). These results are derived
assuming weak conditions on the tail decay of the random
environment.

\bigskip

\noindent\textit{MSC: 60K37;82B44}  

\bigskip

\noindent\textit{Keywords:} Directed polymers in random environment.

\end{abstract}

\section{Introduction}
In this article, we consider a model of directed polymers in random
environment introduced by Huse and Henley in 1985 (\cite{cf:HuHe85}) to
modelize impurity-induced domain-wall roughening in the 2D-Ising
model. This model relates to many physical models of growing random
surfaces including the well known Kardar-Parisi-Zhang equation driven
by gaussian noise (we refer to \cite{cf:KrSp91} for an account on these
models and their relations). In \cite{cf:Zh90}, Zhang proposed to replace
the gaussian noise in the KPZ equation by a noise with power-law tail
to describe fluid flows. Since then, this model has been used to
describe fire fronts, bacterial colonies, etc.... In the field of polymers,
the authors of \cite{cf:HaPa93},\cite{cf:MaZh90} study the random energy landscape of
zero temperature directed polymers in power-law environment distributions.

The first mathematical study of directed polymers at positive
temperature was undertaken by
Imbrie, Spencer in 1988 (\cite{cf:ImSp88}) and carried out by numerous 
authors (\cite{cf:AlZh96},\cite{cf:Bol89},\cite{cf:CaHu02},\cite{cf:CSY03},\cite{cf:Sin},\cite{cf:SoZh96}); for an
overview of the achieved results, we refer to \cite{cf:CSY04}.  In
\cite{cf:CaHu02} and \cite{cf:CSY03}, the authors show, using
martingale techniques,
that the quenched free energy is strictly less than the annealed one
if and only if a localization theorem for the polymer's favorite point
holds. In all these previous mathematical articles, the authors assume
that the environment has exponential moments of all order. When
considering a temperature where the moment generating function of the
environment is infinite, no martingale technique can be used, making
the usual strategy irrelevant. Hence, a natural question in this case
is: what is left from the localization picture? In this
paper, our approach is more general than the martingale approach used
in the above references and we obtain our localization results under much weaker conditions
on the distribution of the environment (including the power tail
distributions studied in \cite{cf:HaPa93},\cite{cf:MaZh90},
exponential distributions...). The case of exponentially distributed
environments is of particular interest in view of the exact results
derived by Johansson in \cite{cf:Jo00} for directed last passage percolation
with i.i.d. exponential variables in dimension $d=1$. Since directed
last passage percolation can be recovered from directed polymers by letting
the temperature go to $0$, one can view the polymer measure as an
interpolation between the directed percolation model and the simple
random walk.

In this article, we go a step further than favorite point localization
and derive localization
results in terms of $\epsilon$-atoms by using bounds on the free
energy. We call  $\epsilon$-atoms, atoms of the polymer measure of
mass at least $\epsilon$. Roughly, we show, under certain assumptions
on the environment, that the whole mass of the polymer measure at "low
temperature" is essentially carried by $\epsilon$-atoms (cf. theorems
\ref{th:lo} and \ref{th:lo"} below). Our method of
proof relies mainly on a simple inequality (cf. lemma
\ref{lem:inf} below) and on an upper bound on greedy lattice animals
established in \cite{cf:Ma02}. Using lemma \ref{lem:inf'}, we also give
a different proof for localization in terms of the
polymer's favorite point if the quenched free energy is strictly
less than the annealed one (cf. theorem \ref{th:lo'} below).

The article is organized as follows: in section 2, we introduce the
model and the definition of $\epsilon$-atoms. In section 3, we state
an existence theorem for the free energy and our localization
theorems. In section 4, we give the proofs.

\section{The model and definition of $\epsilon$-atoms}

\subsection{The model}

The model we consider in this paper consists of a simple random walk
under a random Gibbs measure depending on the temperature. More precisely,

Let $((\omega_{n})_{n \in \N},P)$ denote the simple
random walk starting from $0$ on the $d$-dimensional integer lattice $\Z^d$, defined on
a measurable space $(\Omega, \mathcal{F})$; more precisely, under the measure $P$, $(\omega_{n}-\omega_{n-1})_{n \geq
  1}$ are independent and 
\begin{equation*}
P(\omega_{0}=0)=1, \; \; P(\omega_{n}-\omega_{n-1}=\pm
\delta_{j})=\frac{1}{2d},\; \; j=1, \ldots, d, 
\end{equation*}
where $(\delta_{j})_{1 \leq j \leq d}$ is the j-th vector of the
canonical basis of $\Z^d$.

 The random environment on each lattice site is a sequence
 $\eta=(\eta(n,x))_{(n,x) \in \N \times \Z^d}$ of real valued,
 non-constant and i.i.d. random variables defined on a probability
 space $(H, \mathcal{G}, Q)$. We denote by $F$ the
 common distribution function of the sequence $(\eta(n,x))_{(n,x) \in
 \N \times \Z^d}$. In the whole paper, we will suppose the following:  

\textbf{Assumptions}:
\begin{equation}\label{eq:hyp}
\int_{0}^{\infty}(1-F(x))^{\frac{1}{d+1}}dx < \infty
\end{equation}

and 

\begin{equation}\label{eq:hyp'}
Q(|\eta(n,x)|) < \infty.
\end{equation}
Let $\lambda$ be the logarithmic moment generating function of $\eta(n,x)$:
$$\forall \beta \in \R_{+} \; \; \lambda(\beta)\overset{def.}{=} \ln Q(e^{\beta  \eta(n,x)})\leq \infty.$$

For any $n>0$, we define the (Q-random) polymer measure $\mu_{n}$ on $(\Omega, \mathcal{F})$ by:
$$\mu_{n}(d\omega)=\frac{1}{Z_{n}}\exp(\beta H_{n}(\omega))P(d\omega)$$
where $\beta \in \R_{+}$ is the inverse temperature, 
$$
 H_{n}(\omega)\overset{def.}{=}\sum_{j=1}^{n}\eta(j,\omega_{j}) \qquad
$$
is the hamiltonian and
$$
Z_{n}=P(\exp(\beta H_{n}(\omega))) \qquad
$$
is the partition function.

The above definition shows that the polymer is attracted to sites
where the environment is large and positive and repelled by sites
where the environment is large and negative; as the inverse
temperature $\beta$ increases, the influence of the environment
increases and tends to push the random walk in a few "corridors" where
the environment takes high positive values: we will see in the next
sections quantitative statements of these heuristics.

\subsection{Definition of $\epsilon$-atoms}
The purpose of this article is to study where the polymer measure 
$(\mu_{j-1}(\omega_{j}=x))_{x \in \Z^d}$ is concentrated for large
$j$; under assumptions on the environment $\eta$ and on the inverse
temperature $\beta$ (typically $\beta$ "large"), we show in some sense
that the mass carried by a few points is "significant". To give a
quantitative statement of this phenomenom, we are naturally lead to
introduce the notion of $\epsilon$-atoms. More precisely, let
$\epsilon>0$ be some positive real number; we define
$\mathcal{A}_{j}^{\epsilon,\beta}$ the set of $\epsilon$-atoms to be
the points of $\Z^d$ wich carry a mass of at least $\epsilon$: 

\begin{equation*}
\mathcal{A}_{j}^{\epsilon,\beta}=\lbrace x \in \Z^d:
\mu_{j-1}(\omega_{j}=x)>\epsilon \rbrace.
\end{equation*}

For $\delta<1$, we define the event ${A}_{j}^{\epsilon,\delta,\beta}$
to be the environments for which $\mathcal{A}_{j}^{\epsilon,\beta}$
has a mass of at least $\delta$:

\begin{equation*}
{A}_{j}^{\epsilon,\delta,\beta}=\lbrace\eta:\mu_{j-1}(\omega_{j} \in
\mathcal{A}_{j}^{\epsilon,\beta})  \geq \delta \rbrace,
\end{equation*}

and ${A}_{j}^{\epsilon,\beta}$ to be the environments for which
$\mathcal{A}_{j}^{\epsilon,\beta}$ has at least one element:
\begin{equation*} 
{A}_{j}^{\epsilon,\beta}=\lbrace\eta:\max_{x \in
  \Z^d}\mu_{j-1}(\omega_{j}=x)>\epsilon \rbrace.
\end{equation*}

In terms of $\epsilon$-atoms, we state the following localization
result derived in \cite{cf:CSY03} under the assumption
$\lambda(\beta)<\infty \:(\forall \beta)$ (cf. corollary 2.2 therein):
\begin{equation*}
p(\beta)<\lambda(\beta) \; \Leftrightarrow \; 
\exists \:  \epsilon > 0, \; 
\liminf_{n \to \infty} \frac{1}{n}\sum_{j=1}^{n} \mu_{j-1}(\omega_{j} \in
\mathcal{A}_{j}^{\epsilon,\beta}) \geq \epsilon \qquad  Q-a.s.
\end{equation*}

This equivalence asserts that  the quenched free energy is strictly less
than the annealed one if and only if there exists $\epsilon > 0$ 
such that the mass carried by the $\epsilon$-atoms is for large $n$
(in the sense of C{\'e}saro) bounded below by some positive constant.

We recall that if for all $\beta$ in $\R$ the moment generating
function $\lambda(\beta)$ is finite then for all $\beta$ different
from $0$ the strict inequality $p(\beta)<\lambda(\beta)$ holds for
dimension $d=1$ (theorem 1.1 in \cite{cf:CoVa05}). For dimension $d=2$, this problem is still open.    

Finally, we introduce the following definition 
for $(\mu_{j-1}(\omega_{j} \in .))_{j \geq 1}$:
\begin{definition}
The sequence $(\mu_{j-1}(\omega_{j} \in .))_{j \geq 1}$ is
asymptotically purely atomic (in C{\'e}saro mean) if for all sequence
$(\epsilon_{j})_{j \geq 1}$ tending to $0$ as $j$ goes to infinity the
following convergence holds:  
\begin{equation*}
\frac{1}{n}\sum_{j=1}^{n} \mu_{j-1}(\omega_{j} \in
\mathcal{A}_{j}^{\epsilon_{j},\beta}) \underset{n \to
  \infty}{\longrightarrow} 1   \text{ in Q-Probab.}
\end{equation*}
\end{definition}

Less formally, $(\mu_{j-1}(\omega_{j} \in .))_{j \geq 1}$ is
asymptotically purely atomic if, for large $j$, the polymer measure
concentrates on few atoms.

\section{Results}
\subsection{Existence of the free energy}
First, we establish the existence of the free energy for all $\beta$
in $\R_{+}$. We recall that condition
(\ref{eq:hyp}) implies that $Q[(\eta(n,x)_{+})^{d+1}]<\infty$ but is implied by the existence of some $\epsilon > 0$ such that
$Q[(\eta(n,x)_{+})^{d+1+\epsilon}]< \infty$; in particular, it is much
weaker than the existence of exponential moments. 
We denote by
$\overset{\rightarrow}{\Pi}_{n}$ the oriented paths of the simple
random walk up to time $n$:  
\begin{equation*}
\overset{\rightarrow}{\Pi}_{n}=\lbrace (j,\omega_{j})_{1 \leq j \leq n};
\quad \forall j, \quad |\omega_{j+1}-\omega_{j}|=1 \rbrace
\end{equation*}
We define $\overset{\rightarrow}{N}(n)$ as the maximum of the
environment $\eta$ along the paths of $\overset{\rightarrow}{\Pi}_{n}$: 
\begin{equation*}
\overset{\rightarrow}{N}(n)\overset{def.}{=}\max_{(j,\omega_{j})_{j} \in
  \overset{\rightarrow}{\Pi}_{n}}\sum_{j=1}^{n} \eta(j,\omega_{j}).
\end{equation*}
As a consequence of condition (\ref{eq:hyp}) and proposition 3.4 in \cite{cf:Ma04}, we can define $\alpha$ in
the following way: 
\begin{equation*}
\alpha \overset{def.}{=}\sup_{n \geq 1}Q(\frac{\overset{\rightarrow}{N}(n)}{n}) < \infty.
\end{equation*}
Since $Q(|\eta(n,x)|) < \infty$, by Kingman's subadditive ergodic theorem, we get that: 
\begin{equation*}
\frac{\overset{\rightarrow}{N}(n)}{n}\underset{n \to
  \infty}{\longrightarrow} \alpha \qquad Q-a.s.\text{ and in } L^{1}(Q).
\end{equation*}
For all $\beta$ in $\mathbb{R}_{+}$, the obvious bound $\ln Z_{n}
\leq \beta\overset{\rightarrow}{N}(n)$ and condition (\ref{eq:hyp'})
ensure the existence in $[Q(\eta(n,x)),\alpha]$ of  
\begin{equation*}
p(\beta)=\sup_{n \geq
  1}Q(\frac{\ln Z_{n}}{n}). 
\end{equation*}
To get a strong convergence 
result, we introduce the following condition: 
\begin{equation}\label{eq:hyp"}
\int_{-\infty}^{0}F(x)^{\frac{1}{d+1}}dx < \infty.
\end{equation}
\begin{theorem}\label{th:fe}
The averaged free energy
exists in the following weak sense:   
\begin{equation*}
\frac{Q(\ln Z_{n})}{n} \underset{n \to \infty}{\longrightarrow} p(\beta).
\end{equation*}
We have the following bound on the free energy:
\begin{equation}\label{eq:el}
p(\beta) \leq \alpha \beta \wedge \lambda(\beta).
\end{equation}
If, in addition, the environment satisfies condition (\ref{eq:hyp"}),
one gets the following  stronger result:
\begin{equation*}
\frac{\ln Z_{n}}{n} \underset{n \to \infty}{\longrightarrow} p(\beta)
\qquad Q-a.s.\:and \: in \: L^{1}(Q).
\end{equation*}
\end{theorem}

However, not much is known on the limit $p$: $p$ is convex and
$p(\beta)/\beta\rightarrow
\alpha$ as $\beta \to \infty$. In subsection 3.3, we will tackle the question of the comparaison
of $p$ with its annealed bound $\lambda$.

\subsection{Strong localization in probability}
In this subsection, we fix the inverse temperature $\beta > 0$ and we
suppose that: 
\begin{equation*}
\lambda(\beta)=\infty.
\end{equation*}
Intuitively, when $\lambda(\beta)=\infty$, the environment can take
large values and one expects the polymer measure to concentrate in
those regions of high environment. A quantitative statement of this is
the following theorem:

\begin{theorem}\label{th:lo}
Suppose that $\lambda(\beta)=\infty$. Then, for all $\delta < 1$, there exists $\epsilon(\delta)>0$ such that:
\begin{equation}\label{eq:2}
\liminf_{n \to \infty} Q(\frac{1}{n}\sum_{j=1}^{n} \mu_{j-1}(\omega_{j} \in
\mathcal{A}_{j}^{\epsilon(\delta),\beta})) \geq \delta.
\end{equation}
\end{theorem}


An immediate corollary of the above theorem is the following
convergence result:
\begin{corollary}
The sequence $(\mu_{j-1}(\omega_{j} \in .))_{j \geq 1}$ is
asymptotically purely atomic.
\end{corollary}


\subsection{Almost sure strong localization}
In order to get almost sure localization results, we must suppose that
the environment has non trivial exponential moments. More precisely,  let 
$R= \sup \lbrace \beta \in \R^{+} : \lambda(\beta) < \infty \rbrace$.
In this subsection, we will suppose that $R>0$ (possibly
$R=\infty$). On the interval $]0,R[$, we want to compare $p$ to its
annealed bound $\lambda$, a standard procedure in statistical
physics. Roughly, we have the following conjectured picture for
directed polymers:
\begin{enumerate}
\item
when $p(\beta)=\lambda(\beta)$, $\mu_{n}$ spreads out uniformly.
\item
when $p(\beta)<\lambda(\beta)$, $\mu_{n}$ has macroscopic atoms which
may  concentrate the whole mass.
\end{enumerate}
When $d\geq 3$ and $\beta$ satisfies:
\begin{equation*}
\lambda(2\beta)-2\lambda(\beta)<\ln(1/P(\exists n \geq 1,
  \;\omega_{n}=0)) 
\end{equation*}
(this condition implies $p(\beta)=\lambda(\beta)$),
the situation is well understood: the polymer is
diffusive in the sense that the measure $\mu_{n}(\omega_{n}/\sqrt{n} \in .)$
converges weakly to a gaussian law
(\cite{cf:Bol89},\cite{cf:ImSp88},\cite{cf:SoZh96}) and satisfies a
local limit theorem (\cite{cf:Sin},\cite{cf:Va05}). For case (2), we
refer to theorem \ref{th:lo'} below.

First, we give a preliminary lemma wich states that there is a phase
transition between case (1) and (2) under some assumptions on the environment.

\begin{lemma}\label{lem:dec}
The function $p-\lambda$ is nonincreasing on the interval
$[0,R[$. Suppose that one of the two following conditions is
satisfied:   
\begin{itemize}
\item
$R<\infty$ and $\alpha < \frac{\lambda(R)}{R}$.
\item
$R=\infty$ and defining $L=\text{esssup}(\eta(n,x))$, we have 
$$Q(\eta(n,x)=L)<\overset{\rightarrow}{p}_{c}(d),$$
where $\overset{\rightarrow}{p}_{c}(d)$  denotes the site percolation
threshold for the oriented graph induced on $\mathbb{N}\times
\mathbb{Z}^{d}$ by the simple random walk.
\end{itemize}

Then there exists $\beta_{c}<R$ such that: 
\begin{eqnarray*}
\beta \in [0,\beta_{c}[ & \Rightarrow & p(\beta)=\lambda( \beta). \\ 
\beta \in ]\beta_{c},R[ & \Rightarrow & p(\beta)<\lambda( \beta). \\ 
\end{eqnarray*}
\end{lemma}

\proof
One can adapt the proof of lemma 3.3 in \cite{cf:CoYo06} to prove that
$p-\lambda$ is nonincreasing on the interval $[0,R[$.

If $R<\infty$ and $\alpha < \frac{\lambda(R)}{R}$ then 
\begin{equation*}
\limsup_{\
  \beta \to R}(p(\beta)-\lambda(\beta))<0,
\end{equation*} 
and the existence of $\beta_{c}$ follows.

If $R=\infty$ and $Q(\eta(n,x)=L)<\overset{\rightarrow}{p}_{c}(d)$,
then by standard branching process arguments (cf. theorem 6.1 in
\cite{cf:K84}), one can show that $\alpha<L$. Therefore,  
\begin{equation*}
p(\beta)-\lambda(\beta)\underset{\beta \to \infty}{\sim}\beta(\alpha-L)\underset{\beta \to \infty}{\rightarrow}-\infty 
\end{equation*}
and the existence of $\beta_{c}$ follows.
\qed

\begin{rem}
In lemma \ref{lem:dec}, one can have $\beta_{c}=0$. It is believed
that this is the case in dimension $d=1$ and $d=2$.
\end{rem}

In particular, lemma \ref{lem:dec} gives sufficient conditions for the
existence of $\beta$ in $]0,R[$ such that the strict inequality 
$p(\beta)<\lambda(\beta)$ holds. Now, we state our
first almost sure localization result which generalizes corollary 2.2
in \cite{cf:CSY03}:

\begin{theorem}\label{th:lo'}
Suppose that the environment satisfies condition (\ref{eq:hyp"}). Then for all $\beta$ in $]0,R[$, we have the following implication: 
\begin{equation*}
p(\beta)<\lambda(\beta) \qquad \Rightarrow \qquad \exists 
\: \epsilon > 0, \qquad 
\liminf_{n \to \infty} \frac{1}{n}\sum_{j=1}^{n} \mu_{j-1}(\omega_{j} \in
\mathcal{A}_{j}^{\epsilon,\beta}) \geq \epsilon. \qquad  Q-a.s.
\end{equation*}
\end{theorem}

In the next theorem, we will make the assumption
that $\eta$ "explodes" at $R$:
\begin{equation}\label{eq:ex}
\lambda(R)/R = \infty,
\end{equation}
where wet set $\lambda(R)/R= \text{essup}(\eta(n,x))$ if $R = \infty$ and that
\begin{equation}\label{eq:th}  
\exists \theta > 1, \qquad Q(|\eta(n,x)|^{\theta})<\infty.
\end{equation}

\begin{theorem}\label{th:lo"}
Suppose that the environment satisfies conditions
(\ref{eq:ex}), (\ref{eq:th}). Then for all $\delta<1$, there exists
$\epsilon(\delta)>0$ and $\beta(\delta)$ in $]0,R[$ such that: 
\begin{equation*}
\forall \beta \in [\beta(\delta),R[ \qquad \liminf_{n \to \infty} 
\frac{1}{n}\sum_{j=1}^{n} \mu_{j-1}(\omega_{j} \in
\mathcal{A}_{j}^{\epsilon(\delta),\beta}) \geq \delta \qquad  Q-a.s.
\end{equation*}
\end{theorem}

The above theorem can be seen as a continuity result in view of
theorem \ref{th:lo}.

\section{Proof of theorem \ref{th:fe}}

\noindent\emph{Proof of theorem \ref{th:fe}}
Let $L \in \mathbb{N}^{*}\cup \lbrace \infty \rbrace$. We define
 $Y_{n,L}$ by 
\begin{equation*} 
Y_{n,L}=\frac{1}{n} \ln P(e^{\beta H_{n}^{L}(\omega)}).
\end{equation*}
where 
\begin{equation*}
H_{n}^{L}(\omega)=\sum_{j=1}^{n}\eta(j,\omega_{j})\wedge L
  \vee -L.
\end{equation*}
Similarly, we define
\begin{equation*}
p_{L}(\beta)=\sup_{n \geq 1}Q(Y_{n,L}).
\end{equation*}

With these notations, $Y_{n,\infty}=\frac{1}{n}\ln Z_{n}$ and $p_{\infty}(\beta)=p(\beta)$. It is well
known that the sequence $(Q(\ln Z_{n}))_{n \geq 1}$ is
superadditive and so we have the following limit: 
\begin{equation*} 
\lim_{n \to \infty}Q(Y_{n,\infty})=p(\beta).
\end{equation*}
The obvious bound $\ln Z_{n} \leq \beta \overset{\rightarrow}{N}(n)$
ensures $p(\beta) \leq \alpha \beta$ and an application of Jensen's
inequality to $\ln$ ensures that $p(\beta) \leq \lambda(\beta)$,
giving the first two parts of theorem \ref{th:fe}. From know on, we
suppose the environment satisfies condition (\ref{eq:hyp"}). 

For all $L \in \mathbb{N}^{*}$, it is known (cf. proposition 2.5 in \cite{cf:CSY03}) that 
\begin{equation*}
Y_{n,L}\underset{n \to \infty}{\longrightarrow}p_{L}(\beta) \qquad
 \text{Q-a.s. and in }L^{1}(Q).
\end{equation*}
 One has the following bounds:
\begin{align*} 
\forall n,L \geq 1, \qquad |Y_{n,\infty}-Y_{n,L}| & = |\frac{1}{n}\ln
P(e^{\beta H_{n}(\omega)-H_{n}^{L}(\omega)})| \\
& \leq \frac{\beta}{n} \max_{\omega \in
  \overset{\rightarrow}{\Pi}_{n}}|H_{n}(\omega)-H_{n}^{L}(\omega)| \\ 
& = \frac{\beta}{n}\max_{(j,\omega_{j})_{j} \in
  \overset{\rightarrow}{\Pi}_{n}}\sum_{j=1}^{n} (|\eta|(j,\omega_{j})-L)_{+}.
\end{align*}
Therefore proposition 3.4 in \cite{cf:Ma04} ensures the existence of some
constant $c < \infty$ such that the following estimates hold 
\begin{align*}
\limsup_{n \to \infty}|Y_{n,\infty}-Y_{n,L}| & \leq
c\beta\int_{L}^{\infty}(1-F(x)+F(-x))^{\frac{1}{d+1}}dx \\
& \leq
c\beta\int_{L}^{\infty}(1-F(x))^{\frac{1}{d+1}}dx+c\beta\int_{-\infty}^{-L}F(x)^{\frac{1}{d+1}}dx
\qquad Q-a.s.
\end{align*}
and similarly
\begin{equation*}
\limsup_{n \to \infty}Q(|Y_{n,\infty}-Y_{n,L}|) \leq c\beta\int_{L}^{\infty}(1-F(x))^{\frac{1}{d+1}}dx+c\beta\int_{-\infty}^{-L}F(x)^{\frac{1}{d+1}}dx.
\end{equation*}
Therefore we have
\begin{align*}
|Y_{n,\infty}-Q(Y_{n,\infty})| & \leq
 |Y_{n,\infty}-Y_{n,L}|+|Y_{n,L}-Q(Y_{n,L})|+|Q(Y_{n,L})-Q(Y_{n,\infty})| \\
& \leq
|Y_{n,\infty}-Y_{n,L}|+|Y_{n,L}-Q(Y_{n,L})|+Q(|Y_{n,L}-Y_{n,\infty}|).
\\
\end{align*}
By letting $n \to \infty$, we get 
\begin{equation*}
\limsup_{n \to \infty}|Y_{n,\infty}-p(\beta)| \leq 2c
\beta\int_{L}^{\infty}(1-F(x))^{\frac{1}{d+1}}dx+2c
\beta\int_{-\infty}^{-L}(F(x))^{\frac{1}{d+1}}dx \qquad Q-a.s.
\end{equation*}
By letting $L \to \infty$ above, we conclude 
\begin{equation*}
Y_{n,\infty}\overset{Q-a.s.}{\underset{n \to \infty}{\longrightarrow}}p(\beta).
\end{equation*}
Similarly, we obtain 
\begin{equation*}
Y_{n,\infty}\overset{L^{1}(Q)}{\underset{n \to \infty}{\longrightarrow}}p(\beta).
\end{equation*}
\qed

\section{Proof of theorems \ref{th:lo},\ref{th:lo'},\ref{th:lo"}}

\subsection{Some preliminary lemmas}

We first introduce a few notations we will use in the following two
lemmas. For $n$ a postive integer, we define $\mathcal{P}_{n}$ to be
the standard probability simplex in $\mathbb{R}^{n}$:
\begin{equation*}
\mathcal{P}_{n}=\lbrace(\lambda_{i})_{1 \leq i \leq n} \in
\mathbb{R}_{+}^{n}; \sum_{i=1}^{n}\lambda_{i}=1 \rbrace.
\end{equation*}
For $\epsilon,\delta \in ]0,1[$, we define
\begin{equation*}
\mathcal{P}_{n}^{\epsilon,\delta}=\lbrace (\lambda_{i})_{1 \leq i \leq n} \in
\mathcal{P}_{n}; \sum_{i=1}^{n}\lambda_{i}1_{\lambda_{i}>\epsilon}
\leq \delta \rbrace
\end{equation*}
In this section, $(X_{i})_{i \geq 1}$ will denote  an i.i.d. sequence
of positive random variables on a probability space
$(H,\mathcal{H},P)$ such that:
\begin{equation*} 
E[|\ln X_{1}|] < \infty.
\end{equation*}

\begin{lemma}\label{lem:inf}
Let $\delta \in ]\frac{1}{2},1[$and $\epsilon \in ]0,1-\delta[$ be
such that $\frac{(1-\delta)}{\epsilon}$ is a positive integer.
 We have for all $n \geq  \frac{(1-\delta)}{\epsilon}+1$:
\begin{equation*}
\underset{(\lambda_{i})_{1\! \leq\! i \! \leq \! n} \in 
\mathcal{P}_{n}^{\epsilon,\delta}}{\inf}
E[\ln(\sum_{i=1}^{n}\lambda_{i}X_{i})]=E[\ln(\epsilon \sum_{i=1}^{
  \frac{(1-\delta)}{\epsilon}}X_{i} + \delta X_{\frac{(1-\delta)}{\epsilon}+1})].
\end{equation*}
\end{lemma}
\proof
We can suppose that $X_{1}$ is non constant. We first establish an
auxiliary result we will use intensively in the rest of the proof. Let $k$ be a integer greater than or equal to $2$ and $(\lambda_{i})_{1
  \leq i \leq k}$ an element of $\mathcal{P}_{k}$ such that
$0 < \lambda_{2} \leq \lambda_{1}<1$. 
One can therefore consider the function $\phi:
[0,(1-\lambda_{1}) \wedge \lambda_{2}]
\longrightarrow\mathbb{R}$ defined by:
\begin{equation*} 
\forall \rho \in [0,(1-\lambda_{1})\wedge \lambda_{2}] \qquad \phi(\rho)= E[\ln((\lambda_{1}+\rho)X_{1}+(\lambda_{2}-\rho)X_{1}+\sum_{i=3}^{k}\lambda_{i}X_{i})].
\end{equation*}
One can compute the derivative of $\phi$ and we get $\forall \rho \in ]0,(1-\lambda_{1})\wedge \lambda_{2}]$:
\begin{eqnarray*}
\phi'(\rho) & = & E[\frac{X_{1}-X_{2}}{(\lambda_{1}+\rho)X_{1}+(\lambda_{2}-\rho)X_{2}+\sum_{i=3}^{k}\lambda_{i}X_{i}}]\\
&
= & E[\frac{X_{1}-X_{2}}{(\lambda_{1}+\rho)X_{1}+(\lambda_{2}-\rho)X_{2}+\sum_{i=3}^{k}\lambda_{i}X_{i}}1_{X_{1}>X_{2}}]
\\ & & + E[\frac{X_{1}-X_{2}}{(\lambda_{1}+\rho)X_{1}+(\lambda_{2}-\rho)X_{2}+\sum_{i=3}^{k}\lambda_{i}X_{i}}1_{X_{1}<X_{2}}] \\
& = &
E[\frac{X_{1}-X_{2}}{(\lambda_{1}+\rho)X_{1}+(\lambda_{2}-\rho)X_{2}+\sum_{i=3}^{k}\lambda_{i}X_{i}}1_{X_{1}>X_{2}}]
\\ & & -E[\frac{X_{1}-X_{2}}{(\lambda_{1}+\rho)X_{2}+(\lambda_{2}-\rho)X_{1}+\sum_{i=3}^{k}\lambda_{i}X_{i}}1_{X_{1}>X_{2}}]
\\
& < & 0\\
\end{eqnarray*}
where the last inequality comes from the fact that $x \longrightarrow
\frac{1}{x}$ is decreasing. Therefore, $\phi$ is a decreasing
function and thus one can conclude that $\forall \rho \in ]0,(1-\lambda_{1})\wedge \lambda_{2}]$:
\begin{equation}\label{eq:mon}   
E[\ln((\lambda_{1}+\rho)X_{1}+(\lambda_{2}-\rho)X_{1}+\sum_{i=3}^{k}\lambda_{i}X_{i})]
< E[\ln(\sum_{i=1}^{k}\lambda_{i}X_{i})].
\end{equation}
Let $n \geq  \frac{(1-\delta)}{\epsilon}+1$ be a fixed integer and
consider the application $f: \mathcal{P}_{n}^{\epsilon,\delta} \longrightarrow
\mathbb{R}$ defined by 
\begin{equation*}
\forall (\lambda_{i}) \in \mathcal{P}_{n}^{\epsilon,\delta} \qquad f((\lambda_{i}))=E[\ln(\sum_{i=1}^{n}\lambda_{i}X_{i})].
\end{equation*}
Since $f$ is continuous on the compact set $\mathcal{P}_{n}^{\epsilon,\delta}$,
there exists $(\lambda_{i}^{*}) \in \mathcal{P}_{n}^{\epsilon,\delta}$ such
that:
\begin{equation}\label{eq:inf}
\inf_{(\lambda_{i}) \in \mathcal{P}_{n}^{\epsilon,\delta}} f((\lambda_{i}))=f((\lambda_{i}^{*}))
\end{equation}
Let $p=\# \lbrace i; \lambda_{i}^{*} > 0 \rbrace$. Since $f$ is symmetric,
we can suppose that $\lambda_{1}^{*} \geq \lambda_{2}^{*} \geq \ldots
\lambda_{p}^{*}>0$ and that $\lambda_{i}^{*}=0$ for $i > p$. We introduce the
following set:
\begin{equation*} 
F_{\epsilon}=\lbrace i; \lambda_{i}^{*} > \epsilon \rbrace
\end{equation*}
Let $k=\#F_{\epsilon}$; we have the following identity:
\begin{equation*}
^{c}F_{\epsilon}=[\mid k+1,p \mid]. 
\end{equation*}
If $k \geq 2$, for $\rho>0$ sufficiently small, we have
$(\lambda_{1}^{*}+\rho, \lambda_{2}^{*}-\rho, \lambda_{3}^{*}, \ldots,
\lambda_{p}^{*},0, \ldots, 0) \in A_{n}^{\epsilon,\delta}$ and by inequality
(\ref{eq:mon}), we get 
\begin{equation*}
f(\lambda_{1}^{*}+\rho, \lambda_{2}^{*}-\rho, \lambda_{3}^{*}, \ldots,
\lambda_{p}^{*},0, \ldots, 0) < f((\lambda_{i}^{*})),
\end{equation*}
which contradicts (\ref{eq:inf}). Therefore, $k \leq 1$.  
If $\lambda_{p-1}^{*} < \epsilon$ then for $\rho>0$ sufficiently
small, $(\lambda_{1}^{*}, \ldots,
\lambda_{p-2}^{*},\lambda_{p-1}^{*}+\rho,\lambda_{p}^{*}-\rho,0,
\ldots, 0) \in \mathcal{P}_{n}^{\epsilon,\delta}$ and by inequality
(\ref{eq:mon}) 
we get  
\begin{equation*}
f((\lambda_{1}^{*}, \ldots,
\lambda_{p-2}^{*},\lambda_{p-1}^{*}+\rho,\lambda_{p}^{*}-\rho,0,
\ldots, 0)) < f((\lambda_{i}^{*})),
\end{equation*}
which contradicts (\ref{eq:inf}). Therefore, $(\lambda_{i}^{*})=(\lambda_{1}^{*},
\epsilon, \ldots, \epsilon, \lambda_{p}^{*}, 0, \ldots, 0)$. 
If $\lambda_{1}^{*} < \delta $, then 
for $\rho>0$ sufficiently
small,$(\lambda_{1}^{*}+\rho,\lambda_{2}^{*}-\rho, \lambda_{3}^{*},
\ldots, \lambda_{p}^{*},0,
\ldots, 0) \in \mathcal{P}_{n}^{\epsilon,\delta}$ and by inequality (\ref{eq:mon}) we get  
\begin{equation*}
f((\lambda_{1}^{*}+\rho,\lambda_{2}^{*}-\rho, \lambda_{3}^{*},
\ldots, \lambda_{p}^{*},0,
\ldots, 0)) < f((\lambda_{i}^{*})),
\end{equation*}
which contradicts (\ref{eq:inf}). Thus $\lambda_{1}^{*} = \delta $ and since
$\sum_{i=1}^{p}\lambda_{i}^{*}=1$, we get
$p=1+\frac{1-\delta}{\epsilon}$ and $\lambda_{p}^{*}=\epsilon$. We can 
conclude  
\begin{equation*}
\inf_{(\lambda_{i}) \in \mathcal{P}_{n}^{\epsilon,\delta}} f((\lambda_{i}))=f((\lambda_{i}^{*}))=E[\ln(\epsilon \sum_{i=1}^{
  \frac{(1-\delta)}{\epsilon}}X_{i} + \delta X_{\frac{(1-\delta)}{\epsilon}+1})].\end{equation*}
\qed

\begin{rem}
Under suitable integrability assumptions, the same result holds when
one considers a general concave function instead of $\ln$. 
\end{rem}

In the same spirit than the above lemma, we state the following lemma
without proving it.
\begin{lemma}\label{lem:inf'}
Let $k$ be some positive integer and $\epsilon=\frac{1}{k}$. Then we have for all $n \geq k$:
\begin{equation*}
\underset{\underset{\max(\lambda_{i})\leq \epsilon}{(\lambda_{i})_{1\! \leq\! i \! \leq \! n} \in 
\mathcal{P}_{n}}}{\inf}
E[\ln(\sum_{i=1}^{n}\lambda_{i}X_{i})]=E[\ln(\epsilon \sum_{i=1}^{1/\epsilon}X_{i})].
\end{equation*}

\end{lemma}

Finally, we state the following convergence result:

\begin{lemma}\label{lem:utile}
Let $a,b>0$ be two positive numbers such that $a<b$. We have the
following convergence:
\begin{equation*}
\inf_{\beta \in
  [a,b]}E[\ln(\frac{1}{n}\sum_{i=1}^{n}X_{i}^{\beta})]\underset{n \to
  \infty}{\longrightarrow} \inf_{\beta \in [a,b]}\ln E[X_{1}^{\beta}].
\end{equation*}
\end{lemma}

\proof
The fact that the left hand side is less than or equal to the right
hand side is a consequence of Jensen's inequality.  

Let $L>0$ be such that $-\frac{L}{a}<E[\ln(X_{1})]$. Then for all
$\beta$ in $[a,b]$ we have:
\begin{eqnarray*} 
E[\ln(\frac{1}{n}\sum_{i=1}^{n}X_{i}^{\beta})] & = &
E[\ln(\frac{1}{n}\sum_{i=1}^{n}X_{i}^{\beta})1_{\ln(\frac{1}{n}\sum_{i=1}^{n}X_{i}^{\beta})
  \geq
  -L}] \\ 
& & +E[\ln(\frac{1}{n}\sum_{i=1}^{n}X_{i}^{\beta})1_{\ln(\frac{1}{n}\sum_{i=1}^{n}X_{i}^{\beta}) < -L}] \\ 
& \geq & E[\ln(\frac{1}{n}\sum_{i=1}^{n}(X_{i}\wedge L)^{\beta})1_{\ln(\frac{1}{n}\sum_{i=1}^{n}X_{i}^{\beta})
  \geq
  -L}] \\ 
& & +{\beta}E[\frac{1}{n}\sum_{i=1}^{n}\ln
X_{i}1_{\frac{1}{n}\sum_{i=1}^{n} \ln X_{i} < -\frac{L}{\beta} }]
\\
& \geq & E[\ln(\frac{1}{n}\sum_{i=1}^{n}(X_{i}\wedge L)^{\beta})1_{\ln(\frac{1}{n}\sum_{i=1}^{n}X_{i}^{\beta})
  \geq
  -L}] \\ 
& & +{\beta}E[\ln(X_{1})1_{\frac{1}{n}\sum_{i=1}^{n}\ln(X_{i})<
  -\frac{L}{\beta}}] \\
& \geq & E[\ln(\frac{1}{n}\sum_{i=1}^{n}(X_{i}\wedge L)^{\beta})1_{\ln(\frac{1}{n}\sum_{i=1}^{n}X_{i}^{\beta})
  \geq
  -L}] \\ 
& & -bE[|\ln(X_{1})|1_{\frac{1}{n}\sum_{i=1}^{n}\ln(X_{i})<
  -\frac{L}{a}}]. \\
\end{eqnarray*}
By taking the infimum over all $\beta \in [a,b]$ and using the bounded
convergence theorem, we conclude that:
\begin{equation*}
\liminf_{n \to \infty}\inf_{\beta \in
  [a,b]}E[\ln(\frac{1}{n}\sum_{i=1}^{n}X_{i}^{\beta})] \geq \inf_{\beta \in
  [a,b]}\ln E[(X_{1}\wedge L)^{\beta}].
\end{equation*}
We obtain the result by letting $L \to \infty$ in the above inequality.
\qed

\subsection{Proof of theorem \ref{th:lo}}

Following the notations of lemma \ref{lem:inf}, we consider an
i.i.d. sequence  $(X_{i})_{i \geq 1}$ defined on some probability
space $(H,\mathcal{H},P)$ and such that
$X_{1}\overset{law}{=}e^{\eta(n,x)}$.
Let $\delta < 1$ and
$c(\delta)$ be some integer we will choose at the end of the
proof. Finally, we set $\epsilon=\frac{1-\delta}{c(\delta)}$ (for
notational convenience, we write $\epsilon$ instead of $\epsilon(\delta)$). 

We have the following computation:
\begin{align*}
\frac{Q(\ln Z_{n})}{n}& =
\frac{1}{n}\sum_{j=1}^{n}Q(\ln(\frac{Z_{j}}{Z_{j-1}})) \\
& = \frac{1}{n}\sum_{j=1}^{n}Q(\ln(\sum_{x}\mu_{j-1}(\omega_{j}=x)e^{\beta\eta(j,x)})). \\
& \underset{(Jensen)}{\geq} 
\frac{1}{n}\sum_{j=1}^{n}Q(1_{^{c}A_{j}^{\epsilon,\delta,\beta}}\ln(\sum_{x}\mu_{j-1}(\omega_{j}=x)e^{\beta\eta(j,x)}))+\beta
E[\ln X_{1}]Q(A_{j}^{\epsilon,\delta,\beta}) \\
\end{align*}
Thus, we get the following inequality:
\begin{eqnarray*}
\frac{Q(\ln Z_{n})}{n}-\beta E[\ln X_{1}] & \geq &
\sum_{j=1}^{n}Q(1_{^{c}A_{j}^{\epsilon,\delta,\beta}}\ln(\sum_{x}\mu_{j-1}(\omega_{j}=x)e^{\beta\eta(j,x)}))\\
& & -\beta
E[\ln X_{1}]Q(A_{j}^{\epsilon,\delta,\beta})
\end{eqnarray*}
By applying lemma \ref{lem:inf} to the family $(e^{\beta\eta(j,x)})_{x \in
  \mathbb{Z}^{d}}$ under the conditional measure
  $Q(.|\mathcal{G}_{j-1})$, we get:
\begin{equation*} 
\frac{Q(\ln(Z_{n}))}{n} \geq
(\frac{1}{n}\sum_{j=1}^{n}Q(^{c}A_{j}^{\epsilon,\delta,\beta}))(E[\ln((1-\delta)\sum_{k=1}^{c(\delta)}\frac{1}{c(\delta)}
  X_{k}^{\beta}+\delta 
  X_{c(\delta)+1}^{\beta})]-\beta E[\ln X_{1}]).
\end{equation*}
Therefore, using (\ref{eq:el}) and letting $n$ go to infinity, we get 
\begin{equation*}
\limsup_{n \to
  \infty}(\frac{1}{n}\sum_{j=1}^{n}Q(^{c}A_{j}^{\epsilon,\delta,\beta}))
  \leq \frac{\alpha \beta-\beta E[\ln X_{1}]}{E[\ln((1-\delta)\sum_{k=1}^{c(\delta)}\frac{1}{c(\delta)}
  X_{k}^{\beta}+\delta 
  X_{c(\delta)+1}^{\beta})]-\beta E[\ln X_{1}]}.
\end{equation*}

\begin{equation*}
\limsup_{n \to
  \infty}(\frac{1}{n}\sum_{j=1}^{n}Q(^{c}A_{j}^{\epsilon,\delta,\beta}))
  \leq \frac{\alpha \beta-\beta E[\ln X_{1}]}{E[\ln((1-\delta)\sum_{k=1}^{c(\delta)}\frac{1}{c(\delta)}
  X_{k}^{\beta}+\delta 
  X_{c(\delta)+1}^{\beta})]-\beta E[\ln X_{1}]}.
\end{equation*}
Since $\lambda(\beta) = \infty$, by lemma \ref{lem:utile}, one can choose
$c(\delta)$ such that
\begin{equation*}
\frac{\alpha \beta-\beta E[\ln X_{1}]}{E[\ln((1-\delta)\sum_{k=1}^{c(\delta)}\frac{1}{c(\delta)}
  X_{k}^{\beta}+\delta 
  X_{c(\delta)+1}^{\beta})]-\beta E[\ln X_{1}]} \leq 1- \delta.
\end{equation*}
Since $\mu_{j-1}(\omega_{j} \in
\mathcal{A}_{j}^{\epsilon,\beta})\geq \delta
1_{A_{j}^{\epsilon,\delta,\beta}}$, we get the desired result.

\qed

\subsection{Proof of theorems \ref{th:lo'},\ref{th:lo"}}
Both theorems are based on lemma \ref{lem:inf} or  lemma \ref{lem:inf} and on the law of large numbers
for martingales. Following the notations of lemma \ref{lem:inf}, we consider an
i.i.d. sequence  $(X_{i})_{i \geq 1}$ defined on some probability
space $(H,\mathcal{H},P)$ and such that
$X_{1}\overset{law}{=}e^{\eta(n,x)}$.
We start by proving theorem \ref{th:lo"}. 

\noindent\emph{Proof of theorem \ref{th:lo"}}.
Let $\delta < 1$ and
$c(\delta)$ be some integer we will choose at the end of the
proof. Finally, we set $\epsilon=\frac{1-\delta}{c(\delta)}$ (for
notational convenience, we write $\epsilon$ instead of $\epsilon(\delta)$).
We have the following computation:

\begin{align}
\frac{\ln Z_{n}}{n}& =
\frac{1}{n}\sum_{j=1}^{n}\ln(\frac{Z_{j}}{Z_{j-1}})) \nonumber \\
& =
\frac{1}{n}\sum_{j=1}^{n}\ln(\sum_{x}\mu_{j-1}(\omega_{j}=x)e^{\beta\eta(j,x)}).
\nonumber \\
& = 
\frac{1}{n}\sum_{j=1}^{n}1_{^{c}A_{j}^{\epsilon,\delta,\beta}}\ln(\sum_{x}\mu_{j-1}(\omega_{j}=x)e^{\beta\eta(j,x)})
\nonumber \\
&
+\frac{1}{n}\sum_{j=1}^{n}1_{A_{j}^{\epsilon,\delta,\beta}}\ln(\sum_{x}\mu_{j-1}(\omega_{j}=x)e^{\beta\eta(j,x)}).
\label{eq:pr}\\
\end{align}

Consider the $(\mathcal{G}_{n})$-martingale $M_{n}$ defined by:
\begin{equation*}
M_{n}=\sum_{j=1}^{n}1_{^{c}A_{j}^{\epsilon,\delta,\beta}}(\ln(\sum_{x}\mu_{j-1}(\omega_{j}=x)e^{\beta\eta(j,x)})-Q(\ln(\sum_{x}\mu_{j-1}(\omega_{j}=x)e^{\beta\eta(j,x)})|\mathcal{G}_{j-1})).
\end{equation*}
By definition of $M_{n}$ and by applying lemma \ref{lem:inf} to the
family $(e^{\beta\eta(j,x)})_{x \in \mathbb{Z}^{d}}$ under the
conditional measure $Q(.|\mathcal{G}_{j-1})$, we get: 
\begin{align*}
 & \frac{1}{n}\sum_{j=1}^{n}1_{^{c}A_{j}^{\epsilon,\delta,\beta}}\ln(\sum_{x}\mu_{j-1}(\omega_{j}=x)e^{\beta\eta(j,x)}) \\
& =
M_{n}  
 +\sum_{j=1}^{n}1_{^{c}A_{j}^{\epsilon,\delta,\beta}}Q(\ln(\sum_{x}\mu_{j-1}(\omega_{j}=x)e^{\beta\eta(j,x)})|\mathcal{G}_{j-1}) \\
& \geq M_{n}    
 +(\frac{1}{n}\sum_{j=1}^{n}1_{^{c}A_{j}^{\epsilon,\delta,\beta}})E[\ln((1-\delta)\sum_{k=1}^{c(\delta)}\frac{1}{c(\delta)}
  X_{k}^{\beta}  +\delta 
  X_{c(\delta)+1}^{\beta})]
\end{align*}
Similarly, consider the  $(\mathcal{G}_{n})$-martingale $N_{n}$ defined by:
\begin{equation*}
N_{n}=\sum_{j=1}^{n}1_{A_{j}^{\epsilon,\delta,\beta}}(\ln(\sum_{x}\mu_{j-1}(\omega_{j}=x)e^{\beta\eta(j,x)})-Q(\ln(\sum_{x}\mu_{j-1}(\omega_{j}=x)e^{\beta\eta(j,x)})|\mathcal{G}_{j-1})).
\end{equation*} 
By concavity of ln, we get 
\begin{align*}
&
\frac{1}{n}\sum_{j=1}^{n}1_{A_{j}^{\epsilon,\delta,\beta}}\ln(\sum_{x}\mu_{j-1}(\omega_{j}=x)e^{\beta\eta(j,x)})
 \\
& =  
N_{n}+\sum_{j=1}^{n}1_{A_{j}^{\epsilon,\delta,\beta}}Q(\ln(\sum_{x}\mu_{j-1}(\omega_{j}=x)e^{\beta\eta(j,x)})|\mathcal{G}_{j-1}) \\
& \geq   N_{n}+\beta E[\ln X_{1}](\frac{1}{n}\sum_{j=1}^{n}1_{A_{j}^{\epsilon,\delta,\beta}})
\end{align*}
Plugging the two above inequalities in inequality (\ref{eq:pr}), we get:
\begin{align}
 & \frac{\ln Z_{n}}{n}-\frac{M_{n}}{n}-\frac{N_{n}}{n} -\beta E[\ln X_{1}]  \geq 
\nonumber    \\
& (\frac{1}{n}\sum_{j=1}^{n}1_{^{c}A_{j}^{\epsilon,\delta,\beta}})(E[\ln((1-\delta)\sum_{k=1}^{c(\delta)}\frac{1}{c(\delta)}
  X_{k}^{\beta}+\delta 
  X_{c(\delta)+1}^{\beta})]-\beta E[\ln X_{1}])  \label{eq:pre}
\end{align}
There exists some constant $C>0$ such that for all $j$:
\begin{equation*}  
\beta \sum_{x}\mu_{j-1}(\omega_{j}=x)\eta(j,x)   \leq
\ln(\sum_{x}\mu_{j-1}(\omega_{j}=x)e^{\beta\eta(j,x)}) \leq C
|\sum_{x}\mu_{j-1}(\omega_{j}=x)e^{\beta \eta(j,x)}|^{\frac{1}{\theta}}
\end{equation*}
Thus there exists some constant $C'>0$ such that for all $j$:
\begin{eqnarray*}
Q(|\ln(\sum_{x}\mu_{j-1}(\omega_{j}=x)e^{\beta\eta(j,x)})|^{\theta}) &
\leq & 
C'Q(|\sum_{x}\mu_{j-1}(\omega_{j}=x)e^{\beta\eta(j,x)}|) \\ 
& & +
C'Q(|\sum_{x}\mu_{j-1}(\omega_{j}=x)\eta(j,x)|^{\theta}) \\
& \leq & C'(\lambda(\beta)+Q(|\eta(j,x)|^{\theta})). \\
\end{eqnarray*}
By using theorem 2.19 in \cite{cf:HaHe80}, we conclude that:
\begin{equation*}
\lim_{n \to \infty}\frac{M_{n}}{n}=\lim_{n \to
  \infty}\frac{N_{n}}{n}=0 \qquad Q-a.s.
\end{equation*}
By letting $n$ go to infinity in inequality (\ref{eq:pre}) , we get
\begin{equation*}
\limsup_{n \to \infty}\frac{1}{n}\sum_{j=1}^{n}
1_{^{c}A_{j}^{\epsilon,\delta,\beta}} \leq 
\frac{\alpha \beta-\beta E[\ln X_{1}]}{E[\ln((1-\delta)\sum_{k=1}^{c(\delta)}\frac{1}{c(\delta)}
  X_{k}^{\beta}+\delta 
  X_{c(\delta)+1}^{\beta})]-\beta E[\ln X_{1}]}
\end{equation*}
By using lemma \ref{lem:utile}, one can choose $c(\delta)$
and $\beta(\delta)$ in
$]0,R[$ such that: 
\begin{equation*}
\forall \beta \in [\beta(\delta),R[ \qquad
\frac{\alpha \beta-\beta E[\ln X_{1}]}{E[\ln((1-\delta)\sum_{k=1}^{c(\delta)}\frac{1}{c(\delta)}
  X_{k}^{\beta}+\delta 
  X_{c(\delta)+1}^{\beta})]-\beta E[\ln X_{1}]}\leq 1-\delta,
\end{equation*}
which implies the result since $\mu_{j-1}(\omega_{j} \in
\mathcal{A}_{j}^{\epsilon,\beta})\geq \delta
1_{A_{j}^{\epsilon,\delta,\beta}}$.

\qed

The proof of theorem $\ref{th:lo'}$ follows a similar strategy to the
proof of theorem $\ref{th:lo"}$. Therefore, we only give a sketch of
the proof. 

\noindent\emph{Proof of theorem \ref{th:lo'}}.

Suppose that $\beta$ is such that $p(\beta)<\lambda(\beta)$. Then one
can chose a positive integer $k$ sufficiently large for the following
inequality to hold with $\epsilon=\frac{1}{k}$:
\begin{equation*} 
p(\beta)<E[\ln(\epsilon\sum_{i=1}^{1/\epsilon}X_{i}^{\beta})].
\end{equation*}
By the same strategy than for the proof of theorem \ref{th:lo"} (using
lemma \ref{lem:inf'} instead of lemma \ref{lem:inf}), we
get  

\begin{equation}\label{eq:rec}
\limsup_{n \to \infty}\frac{1}{n}\sum_{j=1}^{n}
1_{^{c}A_{j}^{\epsilon,\beta}} \leq 
\frac{p(\beta)-\beta E[\ln X_{1}]}{E[\ln(\epsilon\sum_{i=1}^{\frac{1}{\epsilon}}X_{i}^{\beta})]-\beta E[\ln X_{1}]}<1.
\end{equation}
This implies easily the desired result.

\qed

\textbf{Acknowledgements}: I would like to thank my Ph.D. supervisor
Francis Comets for his help and suggestions.


\end{document}